\documentclass[a4paper,11pt,reqno]{amsart}

\usepackage{latexsym,dsfont,amssymb,amsmath,amsthm}

\chardef\bslash=`\\ 

\numberwithin{equation}{section}

\newtheorem{theorem}{Theorem}[section]
\newtheorem{corollary}[theorem]{Corollary}
\newtheorem{lemma}[theorem]{Lemma}
\newtheorem{proposition}[theorem]{Proposition}
\theoremstyle{remark}
\newtheorem{remark}[theorem]{Remark}
\newtheorem{example}[theorem]{Example}
\theoremstyle{definition}
\newtheorem{definition}[theorem]{Definition}

\newcommand\bp{\begin{proof}}
\newcommand\ep{\end{proof}}

\newcommand\2{{\frac{1}{2}}}

\newcommand\Dhat{{\hat\Delta}}

\newcommand\D{{\mathcal D}}
\newcommand\E{{\mathcal E}}
\newcommand\F{{\mathcal F}}
\newcommand\RR{{\mathcal R}}
\newcommand\U{{\mathcal U}}

\newcommand\g{{\mathfrak g}}
\newcommand\h{{\mathfrak h}}
\newcommand\kk{{\mathfrak k}}

\newcommand\Clg{\operatorname{Cl}({\mathfrak g})}
\newcommand\Clk{\operatorname{Cl}({\mathfrak k})}
\newcommand\Clv{\operatorname{Cl}(V)}
\newcommand\GL{\operatorname{GL}}
\newcommand\End{\operatorname{End}}
\newcommand\SU{\operatorname{SU}}
\newcommand\Spin{\operatorname{Spin}}

\newcommand{\ad}{\operatorname{ad}}
\newcommand{\add}{{\widetilde\ad}}

\newcommand{\C}{{\mathbb C}}
\newcommand{\R}{{\mathbb R}}
\newcommand\Sp{{\mathbb S}}

\newcommand\eps{\varepsilon}

\newcommand\enu[1]{\smallskip\newline\makebox[5mm][l]{\rm(#1)}}

\begin{document}

\title{The Dirac operator on compact quantum groups}


\author[S. Neshveyev]{Sergey Neshveyev}
\address{Department of Mathematics, University of Oslo,
P.O. Box 1053 Blindern, N-0316 Oslo, Norway.}
\email{sergeyn@math.uio.no}

\author[L. Tuset]{Lars Tuset}
\address{Faculty of Engineering, Oslo University College,
Cort Adelers st. 30, N-0254 Oslo, Norway.}
\email{Lars.Tuset@iu.hio.no}

\thanks{Supported by the Research Council of Norway.}

\dedicatory{Dedicated to the memory of Gerard J. Murphy}

\begin{abstract}
For the $q$-deformation $G_q$, $0<q<1$, of any simply connected
simple compact Lie group~$G$ we construct an equivariant spectral
triple which is an isospectral deformation of that defined by the
Dirac operator $D$ on $G$. Our quantum Dirac operator $D_q$ is a
unitary twist of~$D$ considered as an element of $U\g\otimes\Clg$.
The commutator of $D_q$ with a regular function on $G_q$ consists of
two parts. One is a twist of a classical commutator and so is
automatically bounded. The second is expressed in terms of the
commutator of the associator with an extension of $D$. We show that
in the case of the Drinfeld associator the latter commutator is also
bounded.
\end{abstract}

\date{March 21, 2007}

\maketitle

\bigskip

\section*{Introduction}

The Dirac operator on Minkowski space was introduced in 1928 by P.
Dirac who sought a first order differential operator with square
equal to the Laplacian and found in effect the fundamental
mechanisms governing spin-half particles obeying Fermi
statistics. Its generalization by Atiyah and Singer plays an
essential role in index theory, mathematical physics and
representation theory, and its axiomatization in terms of spectral
triples is at the heart of Connes' non-commutative
geometry~\cite{Co}.

Quantum groups being quantizations of Poisson Lie groups should by
all accounts be non-commutative manifolds, but it has
proved difficult to put them rigorously into Connes' framework. In
the quest for an appropriate Dirac operator on quantum groups and
their homogeneous spaces basically two approaches have been adopted.

One consists of developing $q$-analogues of standard differential
geometric notions, but this poses several immediate problems.
Firstly, it is not clear what a quantum Clifford algebra should be.
The natural suggestion using braidings~\cite{Du1,Fi1,He} seems
reasonable only when the braiding is sufficiently simple, e.g.~a
Hecke symmetry, see for instance the discussion in~\cite{Fi1}.
Secondly, differential calculi are defined in terms of elements of
quantized universal enveloping algebras that do not act as
derivations and yield unbounded commutators with some regular
functions, see e.g.~\cite{Sch}, rendering this approach successful
so far only for irreducible quantum flag manifolds. The case of the
quantum $2$-sphere is considered in~\cite{O,DS1,SW,M} and the
general case is due to Kr\"ahmer~\cite{Kr}, who circumvents the
Clifford algebra problem at the cost of including operators with
unconventional classical limits, see the remark after~\cite[Prop.
2]{Kr}. Also, spectral triples obtained using this approach cannot
be expected to be regular~\cite{NT}.

The other approach was suggested by Connes and Landi~\cite{CoL} and consists of looking
for isospectral deformations of Dirac operators. To write down
explicitly what this means and to handle such operators efficiently
requires a good understanding of Clebsch-Gordan coefficients and
spectral properties of classical Dirac operators. The first result
in this direction was obtained by Chakraborty and Pal~\cite{CP1}, who
constructed a spectral triple on $\SU_q (2)$ which was then studied in
detail by Connes~\cite{Co1}. Their Dirac operator although not
exactly an isospectral deformation of the classical one, is closely
related~\cite{CP2} to such an operator proposed meanwhile by
D\c{a}browski et al~\cite{DLSSV}. Similar results are obtained for
$2$- and $4$-dimensional quantum spheres~\cite{DD1,DD2}. Finally, a
class of operators on quantum $\SU(N)$ and the odd-dimensional
spheres is introduced in~\cite{CP3}, but again it seems difficult to
single out which of these operators have the right classical limit.

The two approaches should be related in the same way as the
universal $R$-matrix is related to the element $t\in\g\otimes\g$
defined by the symmetric invariant form; although we are not going
to discuss this issue in this paper, see Example~\ref{su2} below.

Our construction of the Dirac operator on the $q$-deformation $G_q$
of a group $G$ is inspired by work of Fiore~\cite{Fi0,Fi2}, brought
to our attention by Ulrich Kr\"ahmer. Let $V$ be a finite
dimensional $\g$-module with fixed invariant symmetric form. The
Clifford algebra $\Clv$ is semi-simple and thus has no nontrivial
deformations. For a fixed algebra isomorphism $\varphi\colon U_h\g
\to U\g[[h]]$, Fiore seeks a map $V[[h]]\to\Clv[[h]]$ which
coincides with the embedding $V\to\Clv$ modulo~$h$ and is
$U_h\g$-equivariant. He shows that such a map can be gotten by using
a twist, i.e. an element $\F\in(U\g\otimes U\g)[[h]]$ such that
$(\varphi\otimes\varphi)\Dhat_h=\F\Dhat\varphi(\cdot)\F^{-1}$.

Consider now the Dirac operator on $G$ regarded as an element
$\D$ of $U\g\otimes\Clg$. In the formal deformation setting define the
quantum Dirac operator as the element of $U_h\g\otimes\Clg$ obtained
by conjugating $(\varphi^{-1}\otimes\iota)(\D)$ by
$(\varphi^{-1}\otimes\add)(\F)$, where $\add\colon U\g\to\Clg$ is
induced by the Lie algebra homomorphism $\g\to{\mathfrak{so}}(\g)$.
This can also be done for real parameters $h=2\log q$, but then
instead of $U_q\g$ and $U\g$ one has to consider appropriate
completions.

Although twists exist, their analytical properties are difficult to
study. From a representation theoretical point of view, however, it
is not the twist but the associator and braiding that matter. More
precisely, by a famous result of Drinfeld~\cite{Dr1,Dr2} the
appropriate braided tensor category of $U_h\g$-modules with braiding
defined by the universal $R$-matrix and with trivial associativity
morphisms is equivalent to a category of $U\g[[h]]$-modules with
braiding given by~$e^{ht/2}$ and associativity constraints defined
by the monodromy of the Knizhnik-Zamolodchikov (KZ) equations. Then
choosing a twist essentially means that one fixes such an
equivalence of categories. In the same spirit we show that analytic
properties of our quantum Dirac operator are determined by the
associator rather than the twist, and indeed a form of the
KZ-equations is used crucially to show that we actually get a
spectral triple.

\bigskip

\section{Drinfeld associator}

Let $G$ be a simply connected simple compact Lie group, $\g$ its
complexified Lie algebra. Fix a maximal torus in $G$, and let
$\h\subset\g$ be the corresponding Cartan subalgebra. Choose a
system $\{\alpha_1,\dots,\alpha_n\}$ of simple roots. Let
$(a_{ij})_{1\le i,j\le n}$ be the Cartan matrix of $\g$, and
$d_1,\dots,d_n$ be coprime positive integers such that
$(d_ia_{ij})_{i,j}$ is symmetric. Define as usual a bilinear form on
$\h^*$ by $(\alpha_i,\alpha_j) =d_ia_{ij}$. For each integral
dominant weight $\lambda$ we fix an irreducible unitary
representation $\pi_\lambda\colon G\to B(V_\lambda)$ with highest
weight $\lambda$. Then the group von Neumann algebra~$W^*(G)$ of~$G$ is
the C$^*$-product of the algebras~$B(V_\lambda)$. The algebra
$\U(G)$ of unbounded operators affiliated with $W^*(G)$ is the
algebraic product $\prod_\lambda B(V_\lambda)$. We denote by $\Dhat$
the comultiplication $W^*(G)\to W^*(G)\bar\otimes W^*(G)$. It
extends to a $*$-homomorphism $\U(G)\to \U(G\times
G)=\prod_{\lambda,\mu}B(V_\lambda\otimes V_\mu)$ which we denote by
the same symbol.

For $q\in(0,1)$ denote by $G_q$ the $q$-deformation of $G$. To fix
notation, recall that the algebra~$\U(G_q)$ of unbounded operators
affiliated with the von Neumann algebra $W^*(G_q)$ contains the
algebra $U_q\g$ generated by $X_i$, $Y_i$, $K_i$, $K_i^{-1}$ ($1\le
i\le n$) such that the relations
$$
K_iK_i^{-1}=K_i^{-1}K_i=1,\ \ K_iK_j=K_jK_i,\ \
K_iX_jK_i^{-1}=q_i^{a_{ij}/2}X_j,\ \
K_iY_jK_i^{-1}=q_i^{-a_{ij}/2}Y_j,
$$
$$
X_iY_j-Y_jX_i=\delta_{ij}\frac{K_i^2-K_i^{-2}}{q_i-q_i^{-1}}
$$
as well as the quantum Serre relations are satisfied, where
$q_i=q^{d_i}$. The algebra $U_q\g$ is a Hopf $*$-algebra with
comultiplication $\Dhat_q$ and involution given by
$$
\Dhat_q(K_i)=K_i\otimes K_i,\ \
\Dhat_q(X_i)=X_i\otimes K_i+ K_i^{-1}\otimes X_i,\ \
\Dhat_q(Y_i)=Y_i\otimes K_i+ K_i^{-1}\otimes Y_i,
$$
$$
K_i^*=K_i,\ \ X_i^*=Y_i.
$$

Denote by $\RR\in\U(G_q\times G_q)$ the universal $R$-matrix, see
e.g. \cite[Theorem~8.3.9]{CP} for an explicit formula. It is the
unique element satisfying the following two properties. We have
$$
\Dhat^{op}_q=\RR\Dhat_q(\cdot)\RR^{-1},
$$
and if $\pi_{\lambda,q}$ is a finite dimensional representation with
a highest weight vector $\xi_\lambda$ (that is,
$\pi_{\lambda,q}(X_i)\xi_\lambda =0$ and $\pi_{\lambda,q}(K_i)
\xi_\lambda=q_i^{\lambda(H_i)/2}\xi_\lambda$, where $H_i\in\h$ is
such that $\alpha_j(H_i)=a_{ij}$) and $\pi'_{\mu,q}$ a finite
dimensional representation with a lowest weight vector $\xi'_\mu$
(so $\pi'_{\mu,q}(Y_i)\xi'_\mu=0$ and $\pi'_{\mu,q}(K_i)
\xi'_\mu=q_i^{\mu(H_i)/2}\xi'_\mu$) then
\begin{equation*} \label{eRmat}
(\pi_{\lambda,q}\otimes\pi'_{\mu,q})(\RR)(\xi_\lambda\otimes\xi'_\mu)
=q^{(\lambda,\mu)}\xi_\lambda\otimes\xi'_\mu.
\end{equation*}

Since finite dimensional representations of $G_q$ are again
classified by integral dominant weights, we have a canonical
identification of the centers of $W^*(G_q)$ and $W^*(G)$. It extends
to a $*$-isomorphism $W^*(G_q)\cong W^*(G)$, and therefore
$\U(G_q)\cong\U(G)$. Such an isomorphism does not respect
comultiplications, and to compare them we recall the notion of the
Drinfeld associator, see e.g.~\cite{EFK} and \cite{Ka} for more
details.

\smallskip

Let $A$ and $B$ be operators on a finite dimensional vector space
$V$. Put
$$
\hbar=\frac{\log q}{\pi i}.
$$
Consider the differential equation
$$
G'(x)=\hbar\left(\frac{A}{x}+\frac{B}{x-1}\right)G(x),
$$
where $G\colon(0,1)\to\End(V)$. Assume that neither $A$ nor $B$ has
eigenvalues which differ by a nonzero integral multiple of
$\frac{1}{\hbar}$. Then there exist unique solutions $G_0$ and $G_1$
such that the functions $G_0(x)x^{-\hbar A}$ and $G_1(1-x)x^{-\hbar
B}$ extend to holomorphic functions in the unit disc with value $1$
at $x=0$. These solutions are in fact $\GL(V)$-valued, hence there
exists $\Phi(A,B)\in \GL(V)$ such that
$$
G_0(x)=G_1(x)\Phi(A,B)\ \ \hbox{for all}\ \ x\in(0,1).
$$

We will be interested only in the case when $V$ is a Hilbert space
and the operators $A$ and~$B$ are self-adjoint. Then the assumptions
on the spectra are automatically satisfied. For $a\in(0,1)$ let
$G_a$ be the unique solution such that $G_a(a)=1$. Note that
$G_a(x)$ is unitary. Indeed, $\hbar(\frac{A}{x}+\frac{B}{x-1})$ is
skew-adjoint, so $G_a$ is an integral curve of a time-dependent
vector field on the unitary group. By uniqueness of solutions we
have $G_a(x)=G_0(x)G_0(a)^{-1}$, so
$$
a^{-\hbar B}G_a(1-a)a^{\hbar A} =a^{-\hbar
B}G_0(1-a)G_0(a)^{-1}a^{\hbar A} =a^{-\hbar
B}G_1(1-a)\Phi(A,B)G_0(a)^{-1}a^{\hbar A}.
$$
Since $a^{-\hbar B}$ is unitary for any $a\in(0,1)$, the operators
$$
a^{-\hbar B}G_1(1-a)=a^{-\hbar B}(G_1(1-a)a^{-\hbar B})a^{\hbar B}
$$
converge to $1$ as $a\to0^+$. Similarly $G_0(a)^{-1}a^{\hbar A}\to1$
as $a\to0^+$. It follows that
\begin{equation} \label{eKZ}
\Phi(A,B)=\lim_{a\to 0^+}a^{-\hbar B}G_a(1-a)a^{\hbar A}.
\end{equation}
This expression makes it in particular obvious that $\Phi(A,B)$ is
unitary.

\smallskip

Consider the rescaling $(\cdot,\cdot)$ of the Killing form on $\g$
such that its restriction to $\h$ is  the one induced by the
symmetric form on $\h^*$ defined above. Let $\{x_k\}_k$ be a basis
in the real Lie algebra of $G$ such that $(x_k,x_l)=-\delta_{kl}$.
Put
$$
t=-\sum_k x_k\otimes x_k\in \g\otimes\g\subset\U(G\times G).
$$
This element is self-adjoint, e.g. because $x_k$ lie in the real Lie
algebra of $G$ and so $x_k^*=-x_k$. The Drinfeld associator is
defined by
$$
\Phi_{KZ}=\Phi(t_{12},t_{23}).
$$
More precisely, it is the unique unitary element in
$W^*(G)\bar\otimes W^*(G)\bar\otimes W^*(G)$
such that for any finite dimensional representations $\pi_i\colon
G\to B(V_i)$, $i=1,2,3$, we have
$$
(\pi_1\otimes\pi_2\otimes\pi_3)(\Phi_{KZ})=
\Phi((\pi_1\otimes\pi_2\otimes\pi_3)(t\otimes1),
(\pi_1\otimes\pi_2\otimes\pi_3)(1\otimes t)).
$$

\medskip

The following variant of a famous result of Drinfeld will play a
central role in the paper.

\begin{theorem} \label{Dr}
There exist a $*$-isomorphism $\varphi\colon W^*(G_q)\to W^*(G)$
extending the canonical identification of the centers and a unitary
$\F\in W^*(G)\bar\otimes W^*(G)$ such that \enu{i}
$(\varphi\otimes\varphi)\Dhat_q=\F\Dhat\varphi(\cdot)\F^{-1}$;
\enu{ii} $(\hat\eps\otimes\iota)(\F)=(\iota\otimes\hat\eps)(\F)=1$, where
$\hat\eps$ is the trivial representation of $G$; \enu{iii}
$(\varphi\otimes\varphi)(\RR)=\F_{21}q^t\F^{-1}$; \enu{iv} the
associator $\Phi=(\iota\otimes\Dhat)(\F^{-1})
(1\otimes\F^{-1})(\F\otimes1)(\Dhat\otimes\iota)(\F)$ coincides with
the Drinfeld associator~$\Phi_{KZ}$.
\end{theorem}

If $\F\in\U(G\times G)$ is an element satisfying (i) for some
isomorphism $\varphi\colon\U(G_q)\to \U(G)$ extending the
identification of the centers, then we say that $\F$ is a twist. If
in addition $\varphi$ is a $*$-homomorphism and $\F$ is unitary, we
say that $\F$ is a unitary twist. If all four conditions (i)-(iv)
are satisfied, we talk about Drinfeld twists and unitary Drinfeld
twists.

\bp[Proof of Theorem~\ref{Dr}] The existence of a Drinfeld twist in
the formal deformation setting is due to Drinfeld~\cite{Dr1,Dr2}.
There it is proved by inductive cohomological arguments, so the
twist makes sense only as a formal power series and a priori cannot
be specialized to a complex deformation parameter. The result
implies equivalence of certain braided tensor categories. A
constructive proof of this equivalence was later given by Kazhdan
and Lusztig~\cite{KL1,KL2}, the advantage being that the
specialization makes sense (for nonzero complex parameters different
from nontrivial roots of unity). It also implies the existence of a
Drinfeld twist. The construction of Kazhdan and Lusztig was further
clarified and extended by Etingof and Kazhdan~\cite{EK1,EK2,EK3}.

Therefore there exist an isomorphism $\varphi\colon \U(G_q)\to\U(G)$
and a Drinfeld twist $\F\in\U(G\times G)$, and the only additional
claim we make is that one can choose $\varphi$ to be $*$-preserving
and $\F$ unitary.

\smallskip

Let us show first that $\varphi$ can be assumed to be
$*$-preserving. Since every homomorphism of full matrix algebras is
equivalent to a $*$-homomorphism, there exists an invertible element
$u\in\U(G)$ such that the homomorphism
$$
\varphi_u=u\varphi(\cdot)u^{-1}
$$
is $*$-preserving. We may assume $\hat\eps(u)=1$. Then
$$
\F_u=(u\otimes u)\F\Dhat(u^{-1})
$$
is a Drinfeld twist for $\varphi_u$. Indeed, the conditions (i) and
(ii) are obviously satisfied. To show~(iii) recall that $t$ is
$\g$-invariant,  i.e. $[t,\Dhat(x)]=0$ for any $x$. In particular,
$\Dhat(u)$ commutes with $t$ and recalling that $\Dhat^{op}=\Dhat$
we get
\begin{align*}
(\varphi_u\otimes\varphi_u)(\RR)
&=(u\otimes u)\F_{21}q^t\F^{-1}(u^{-1}\otimes u^{-1})\\
&=(u\otimes u)\F_{21}\Dhat(u^{-1}) q^t\Dhat(u)\F^{-1}(u^{-1}\otimes
u^{-1}) =(\F_u)_{21}q^t\F^{-1}_u.
\end{align*}
Finally, a direct computation shows that the new associator
$$
\Phi_u=(\iota\otimes\Dhat)(\F^{-1}_u)
(1\otimes\F^{-1}_u)(\F_u\otimes1)(\Dhat\otimes\iota)(\F_u)
$$
equals $\Dhat^{(2)}(u)\Phi\Dhat^{(2)}(u^{-1})$. It remains to recall
that $\Phi$ is $\g$-invariant, since $\Dhat_q$ is coassociative.
This is also clear by definition of $\Phi_{KZ}$, as $t$ is
$\g$-invariant and hence $\Phi(t_{12},t_{23})$ is $\g$-invariant.

\smallskip

Assuming now that $\varphi\colon\U(G_q)\to\U(G)$ is a
$*$-isomorphism and $\E$ a Drinfeld twist, we assert that the
unitary $\F$ in the polar decomposition $\E=\F|\E|$ is a unitary
Drinfeld twist for $\varphi$. Indeed, since $\Dhat_q$, $\Dhat$ and
$\varphi$ are $*$-homomorphisms, condition (i) on $\E$ implies that
$$
\E\Dhat\varphi(\cdot)\E^{-1}=(\E^{-1})^*\Dhat\varphi(\cdot)\E^*,
$$
that is, $\E^*\E$ is $\g$-invariant. It follows that $|\E|$ is also
$\g$-invariant. Hence
$$
\E\Dhat\varphi(\cdot)\E^{-1}=\F|\E|\Dhat\varphi(\cdot)|\E|^{-1}\F^{-1}
=\F\Dhat\varphi(\cdot)\F^{-1},
$$
so condition (i) for $\F$ is satisfied. Condition (ii) is also
obviously satisfied. Turning to (iii) recall that the $R$-matrix has
the property $\RR^*=\RR_{21}$. So applying the $*$-operation and
then the flip to the identity
$(\varphi\otimes\varphi)(\RR)=\E_{21}q^t\E^{-1}$ we get
$(\varphi\otimes\varphi)(\RR)=(\E^{-1})^*_{21}q^t\E^*$. Therefore
$$
(\E^*\E)_{21}q^t=q^t\E^*\E
$$
and hence $|\E|_{21}q^t=q^t|\E|$. It follows that
$$
(\varphi\otimes\varphi)(\RR)=\E_{21}q^t\E^{-1}
=\F_{21}|\E|_{21}q^t|\E|^{-1}\F^{-1}=\F_{21}q^t\F^{-1}.
$$
It remains to check (iv). Consider the new associator
$$
\Phi_0=(\iota\otimes\Dhat)(\F^{-1})
(1\otimes\F^{-1})(\F\otimes1)(\Dhat\otimes\iota)(\F).
$$
We have to show that $\Phi_0=\Phi$. Since $|\E|$ is $\g$-invariant,
one easily checks that
\begin{equation} \label{eAss}
\Phi=(\iota\otimes\Dhat)(|\E|^{-1})
(1\otimes|\E|^{-1})\Phi_0(|\E|\otimes1)(\Dhat\otimes\iota)(|\E|).
\end{equation}
Since $\Phi_0$ is defined by the unitary element $\F$, it is itself
unitary. Since $\Phi=\Phi_{KZ}$ is also unitary, taking the inverses
in the above identity and then applying the $*$-operation we get
$$
\Phi=(\iota\otimes\Dhat)(|\E|)
(1\otimes|\E|)\Phi_0(|\E|^{-1}\otimes1)(\Dhat\otimes\iota)(|\E|^{-1}).
$$
Therefore
$$
(\iota\otimes\Dhat)(|\E|^{-1})
(1\otimes|\E|^{-1})\Phi_0(|\E|\otimes1)(\Dhat\otimes\iota)(|\E|)
=(\iota\otimes\Dhat)(|\E|)
(1\otimes|\E|)\Phi_0(|\E|^{-1}\otimes1)(\Dhat\otimes\iota)(|\E|^{-1}).
$$
Since $|\E|$ is $\g$-invariant, the positive operators
$(\iota\otimes\Dhat)(|\E|)$ and $1\otimes|\E|$, as well as
$|\E|\otimes1$ and $(\Dhat\otimes\iota)(|\E|)$, commute. So we can
write
$$
\Phi_0((|\E|\otimes1)(\Dhat\otimes\iota)(|\E|))^2
=((1\otimes|\E|)(\iota\otimes\Dhat)(|\E|))^2\Phi_0.
$$
Consequently
$$
\Phi_0(|\E|\otimes1)(\Dhat\otimes\iota)(|\E|)
=(1\otimes|\E|)(\iota\otimes\Dhat)(|\E|)\Phi_0,
$$
and returning to (\ref{eAss}) we get $\Phi=\Phi_0$. \ep

\begin{remark}
For any $*$-isomorphism $\varphi\colon W^*(G_q)\to W^*(G)$ extending
the identification of the centers, the existence of a unitary twist
satisfying also condition (ii) follows immediately from the fact
that the fusion rules for $G$ and $G_q$ are the same. Then one can
modify the twist to satisfy condition (iii) using the symmetrization
procedure of Drinfeld together with the identity
\begin{equation} \label{eRC0}
\RR^*\RR=\Dhat_q(q^{C_q})(q^{-C_q}\otimes q^{-C_q}),
\end{equation}
where $C_q=\varphi^{-1}(C)$ and $C=-\sum_kx_k^2$ is the Casimir, and
the identity
\begin{equation} \label{eCas1}
t=\2(\Dhat(C)-1\otimes C-C\otimes 1).
\end{equation}
Note in passing that the above two identities imply that
\begin{equation} \label{eRC}
(\varphi\otimes\varphi)(\RR^*\RR)=\F q^{2t}\F^*
\end{equation}
holds for any unitary twist $\F$.

Therefore the key condition is number (iv). In the formal
deformation setting Drinfeld proves a stronger result, so it makes
sense to ask the following question. Given a unitary twist~$\F$
satisfying conditions (i)-(iii), does there exist a $\g$-invariant
unitary~$\E$ such that $\E_{21}=\E$ and the associator defined by
$\F\E$ coincides with the Drinfeld associator?
\end{remark}

\bigskip

\section{Dirac operator}

Recall that the Dirac operator $D$ on a spin manifold $M$ is defined
as the composition
$$
\Gamma(S)\xrightarrow{\tilde\nabla}\Gamma(T^*M\otimes
S)\xrightarrow{\sim} \Gamma(TM\otimes S)\xrightarrow{c}\Gamma(S)
$$
of the Clifford action $c$ on the spin bundle $S$ with the spin
connection $\tilde\nabla$, using the metric to identify tangent and
cotangent bundles. Thus with respect to an orthonormal local
frame~$\{e_i\}_i$ the Dirac operator is given by
$D=\sum_ic(e_i)\tilde\nabla_{e_i}$.

Turning to the Dirac operator on $G$, trivialize the complexified
tangent bundle $TG$ by left translations and identify
$\Gamma(TG)$ with $C^\infty(G)\otimes\g$. Define a Riemannian metric
on $G$ using the form on $\g$ introduced earlier multiplied by $-1$.
The Levi-Civita connection is given by
$$
\nabla_{f\otimes x}=f\partial(x)\otimes 1+\2 f\otimes\ad(x),
$$
where $\partial$ is the representation of $U\g$ by left-invariant
differential operators.

Denote by $\Clg$ the complex Clifford algebra of $\g$ and by
$\gamma\colon\g\to\Clg$ the canonical embedding, so $\Clg$ is
generated by $\gamma(x)$, $x\in\g$, and $\gamma(x)^2=(x,x)1$. We
regard $\Clg$ as a $*$-algebra by requiring the map $\gamma$ to be
$*$-preserving. The spin group $\Spin(\g)$ is the connected Lie
subgroup of the group of invertible elements of $\Clg$ with real Lie
algebra spanned by the elements $\gamma(x_k)\gamma(x_l)$, $k\ne l$.
It acts on $\Clg$ by inner automorphisms. The adjoint action of $G$
on $\g$ extends to an action of $G$ on~$\Clg$ which lifts to a
homomorphism $G\to\Spin(\g)$. At the Lie algebra level it is given
by
$$
\g\ni x\mapsto\add(x) :=\frac{1}{4}\sum_k\gamma(x_k)\gamma([x,x_k]).
$$
We denote by the same symbol $\add$ the corresponding homomorphism
$\U(G)\to\Clg$. Note that by definition the map $\gamma$ is
equivariant, so $\gamma([x,y]) =[\add(x),\gamma(y)]$ for $x,y\in\g$.

Fix a spin module, that is, an irreducible $*$-representation
$s\colon\Clg\to B(\Sp)$. Recall that if $\g$ is even dimensional
then $s$ is unique up to equivalence and faithful, and there are two
possibilities for~$s$ in the odd dimensional case. Identifying the
smooth sections of the spin bundle $S=G\times\Sp$ with
$C^\infty(G)\otimes\Sp$, the spin connection is
$$
\tilde\nabla_{f\otimes x}=f\partial(x) \otimes s(1)+\2 f\otimes
s\,\add(x),
$$
The Clifford action of $1\otimes x$ is given by $1\otimes
s\gamma(x)$. Hence using the orthonormal global frame
$\{e_k=1\otimes x_k\}_k$, we see that the Dirac operator $D\colon
C^\infty(G)\otimes\Sp\to C^\infty(G)\otimes\Sp$ is given by
$$
D=\sum_k\left(\partial(x_k)\otimes s\gamma(x_k) +\2\otimes
s(\gamma(x_k) \add(x_k))\right).
$$
This can be written as $D=(\partial\otimes s)(\D)$, where
$$
\D=\sum_k\left(x_k\otimes\gamma(x_k)+\2
\otimes\gamma(x_k)\add(x_k)\right)
$$
is an element of the non-commutative Weil algebra $U\g\otimes\Clg$,
see \cite{AM}.

\begin{remark} \label{Sleb}
One can use other connections than the Levi-Civita one to define a
Dirac operator by varying the coefficient $\2$ in the above
expressions~\cite{Sl}. Taking $0$ one gets an operator corresponding
to the reductive connection, and taking $\frac{1}{3}$ one gets
Kostant's cubic Dirac operator~\cite{Ko}.
\end{remark}

Now fix a unitary twist $\F$ corresponding to a $*$-isomorphism
$\varphi\colon W^*(G_q)\to W^*(G)$. Define the universal quantum
Dirac operator $\D_q\in\U(G_q)\otimes\Clg$ by
$$
\D_q=(\varphi^{-1}\otimes\iota)((\iota\otimes\add)(\F)
\D(\iota\otimes\add)(\F^*)).
$$
Denote by $\C[G_q]$ the linear span of matrix coefficients of finite
dimensional representations of $G_q$. It is a Hopf $*$-algebra with
comultiplication $\Delta_q$, and $\U(G_q)$ is its dual space. Let
$(L^2(G_q),\pi_{r,q},\xi_q)$ be the GNS-triple defined by the Haar
state on $\C[G_q]$. The left $\hat\pi_{r,q}$ and right~$\partial_q$
regular representations of $W^*(G_q)$ on $L^2(G_q)$ are defined by
$$
\hat\pi_{r,q}(\omega)\pi_{r,q}(a)\xi_q=(\omega
S^{-1}\otimes\pi_{r,q})\Delta_q(a)\xi_q,
$$
where $S$ is the antipode on $\C[G_q]$, and
\begin{equation} \label{eRegL}
\partial_q(\omega)\pi_{r,q}(a)\xi_q
=(\pi_{r,q}\otimes\omega)\Delta_q(a)\xi_q
=a_{(1)}(\omega)\pi_{r,q}(a_{(0)})\xi_q.
\end{equation}

\begin{definition}
The quantum Dirac operator $D_q$ is the unbounded operator on
$L^2(G_q)\otimes\Sp$ defined by
$$
D_q=(\partial_q\otimes s)(\D_q).
$$
\end{definition}

\begin{remark} \label{uniq}
The element $\D_q$ depends a priori on the choice of $\varphi$ and
$\F$. Is it true that $\D$ commutes with all $\g$-invariant elements
in $(\iota\otimes\add) (U\g\otimes U\g)$? This is the case for
$\SU(2)$, see Example~\ref{su2} below, but this case is special
since then $s\gamma$ coincides with $s\,\add$ up to a scalar factor.
If the answer is yes in general then $\D_q$ does not depend on $\F$
for fixed $\varphi$. On the other hand, the dependence on $\varphi$
is very mild. Namely, if we replace $\varphi$ by $\varphi'$ then
there exists a unitary $v\in W^*(G)$ such that
$\varphi'=v\varphi(\cdot)v^*$ and $\F_v=(v\otimes v)\F\Dhat(v^*)$ is
a unitary twist for $\varphi'$. Since $\D$ commutes with the image
of $(\iota\otimes\add)\Dhat$, for the element $\D_q'$ defined by
$\varphi'$ and $\F_v$ we get
$\D'_q=(1\otimes\add(v))\D_q(1\otimes\add(v^*))$.
\end{remark}

The quantum group $G_q$ acts on itself from the left and from the
right, and the operator~$D_q$ is equivariant with respect to these
two actions. More formally, we have two coactions of~$\C[G_q]$ on
itself. They can be implemented by the representations
$\hat\pi_{r,q}(\cdot)\otimes1$ and $\partial_q\times
s\,\add_q=(\partial_q\otimes s\,\add_q)\Dhat_q$ of $W^*(G_q)$, where
$\add_q=\add\,\varphi\colon W^*(G_q)\to\Clg$. Then we have the
following.

\begin{proposition}
The universal quantum Dirac operator $\D_q\in\U(G_q)\otimes\Clg$
commutes with all elements of the form
$(\iota\otimes\add_q)\Dhat_q(x)$, where $x\in W^*(G_q)$. In
particular, the quantum Dirac operator~$D_q$ commutes with all
operators of the form $\hat\pi_{r,q}(x)\otimes1$ and
$(\partial_q\times s\,\add_q)(x)$.
\end{proposition}

\bp Recall that as $t=-\sum_kx_k\otimes x_k$ is $\g$-invariant and
the map $\gamma$ is equivariant, the element
$(\iota\otimes\gamma)(t)$ commutes with any element of the form
$(\iota\otimes\add)\Dhat(x)$, and similarly $\sum_k\gamma(x_k)
\add(x_k)$ commutes with any element in the image of $\add$. So $\D$
commutes with any element of the form $(\iota\otimes\add)\Dhat(x)$.
Thus $(\varphi\otimes\iota)(\D_q)$ commutes with any element of the
form
$$
(\iota\otimes\add)(\F\Dhat\varphi(x)\F^*)
=(\varphi\otimes\add\,\varphi)\Dhat_q(x)
=(\varphi\otimes\add_q)\Dhat_q(x),\ \ x\in\U(G_q),
$$
so $\D_q$ commutes with $(\iota\otimes\add_q)\Dhat_q(x)$.

By applying $\partial_q\otimes s$ we see that $D_q$ commutes with
$(\partial_q\times s\,\add_q)(x)$. Finally $D_q$ commutes with
$\hat\pi_{r,q}(x)\otimes1$ simply because $\hat\pi_{r,q}(x)$
commutes with $\partial_q(y)$ for all $y$. \ep

Next note that by definition the operator $D_q$ is unitarily
equivalent to $D$. In particular, $D_q$ is self-adjoint and its
spectrum is the same as that of $D$. Recall that one can compute the
squares of the eigenvalues of $D$ by using the Weitzenb\"ock
formula:
$$
\D^2=\2(\iota\otimes\add)\Dhat(C)+\2 C\otimes 1
+\frac{1}{4}\otimes\add(C).
$$
Recall also that $\add(C)=3\|\rho\|^2$, where $\rho$ is half the sum
of the positive roots, which can be seen using the well-known result
of Kostant~\cite{Ko0} that the representation $s\,\add$ is
equivalent to several copies of the irreducible representation with
highest weight $\rho$, and that the image of~$C$ under an
irreducible representation with highest weight $\lambda$ is the
scalar $\|\lambda+\rho\|^2-\|\rho\|^2$. Therefore
$$
\D^2=\2(\iota\otimes\add)\Dhat(C)+\2 C\otimes 1
+\frac{3}{4}\|\rho\|^2.
$$
For $\D_q$ this can be reformulated as follows.

\begin{proposition} \label{Weitz}
We have
$$
\D^2_q=\2(\iota\otimes\add_q)\Dhat_q(C_q)+\2 C_q\otimes 1
+\frac{3}{4}\|\rho\|^2.
$$
It follows that
$$
q^{2\D_q^2}=(q^{2C_q+\frac{9}{2}\|\rho\|^2}\otimes1)
(\iota\otimes\add_q) (\RR^*\RR).
$$
\end{proposition}

\bp The first identity follows immediately from definitions and the
Wei\-tzen\-b\"ock formula. The second follows from (\ref{eRC0}) and
the equality $\add_q(C_q)=3\|\rho\|^2$. \ep

The proposition shows that $\D_q^2$ does not depend on the choice of
$\F$. Moreover, one can get an explicit formula for $q^{2\D_q^2}$ in
terms of the generators of $U_q\g$ (recall that $C_q$ can also be
expressed in terms of the $R$-matrix by $q^{-2C_q}=\hat
m_q(\iota\otimes\hat S_q)(\RR^*\RR)$, where $\hat m_q$ and $\hat
S_q$ are the multiplication and the antipode on $\U(G_q)$).

As in the classical case, the Weitzenb\"ock formula allows one to
compute the spectral subspaces of~$D_q^2$. Namely, let $\tilde
V_{\lambda,q}\subset L^2(G_q)$ be the linear span of the matrix
coefficients of an irreducible representation with highest weight
$\lambda$. Then $\tilde V_{\lambda,q}\otimes\Sp$ is
$(\partial_q\times s\,\add_q)(W^*(G_q))$-invariant, and if $V\subset
\tilde V_{\lambda,q}\otimes\Sp$ is an irreducible submodule with
highest weight $\mu$, then $D^2_q$ acts on $V$ as the scalar
$$
\2\|\mu+\rho\|^2+\2\|\lambda+\rho\|^2-\frac{1}{4}\|\rho\|^2.
$$
It is worth recalling that if we use the reductive connection
instead of the Levi-Civita one, then using a similar result we can
compute the spectrum of $D_q$ completely, as then for any
eigenvalue~$\beta$ the number $-\beta$ is again an eigenvalue with
the same multiplicity~\cite{Fe}.

\begin{example} \label{su2}
Consider the simplest case $G=\SU(2)$. Then $\Clg$ can be identified
with the algebra $B(V_\2)\oplus B(V_\2)$ in such a way that
$\add=\pi_\2\oplus\pi_\2$ and
$\gamma(x)=\sqrt{2}(\pi_\2(x),-\pi_\2(x))$, $x\in\g$. Choose
$s\colon B(V_\2)\oplus B(V_\2) \to B(V_\2)$ to be the projection on
the first factor. We then see that the map $(\sqrt{2})^{-1}s\gamma$
coincides with the restriction of $s\,\add$ to $\g$, in particular,
it extends to a homomorphism $\U(G)\to B(V_\2)$. This implies that
the element $\D$ commutes with any element of the form
$(\iota\otimes\add)(\E)$, where $\E$ is a $\g$-invariant unitary,
since $t$ commutes with any such unitary by virtue of (\ref{eCas1}).
Hence for any fixed $\varphi$ the operator $\D_q$ is independent of
the twist $\F$. Therefore by Remark~\ref{uniq} we conclude that
$\D_q$ is unique up to the inner automorphism of
$\U(G_q)\otimes\Clg$ defined by a unitary of the form
$1\otimes\add(u)$, $u\in W^*(G_q)$.

Since in our case $\sum_k\add(x_k^2) =-\add(C)=-\frac{3}{2}$, we
have
$$
(\iota\otimes s)(\D)=-\sqrt{2}(\iota\otimes s\,\add)(t)
-\frac{3\sqrt{2}}{4},
$$
whence
$$
q^{-\sqrt{2}(\iota\otimes s)(\D)} =q^{\frac{3}{2}}(\iota\otimes
s\,\add)(q^{2t}).
$$
Recall that by (\ref{eRC}) we have
$(\varphi\otimes\varphi)(\RR^*\RR)=\F q^{2t}\F^*$. It follows that
$$
q^{-\sqrt{2}D_q} =q^{\frac{3}{2}}(\partial_q\otimes
s\,\add_q)(\RR^*\RR).
$$

To get an explicit expression for $q^{-\sqrt{2}D_q}$, consider the
standard generators $e=X_1$, $f=Y_1$, $k=K_1$ of $U_q\g$. The
representation $s\,\add_q$ is an irreducible representation of spin
$\2$, so with an appropriate choice of basis we have
$$
e\mapsto
\begin{pmatrix}
0 & 1\\ 0 & 0
\end{pmatrix},\ \
f\mapsto
\begin{pmatrix}
0 & 0\\ 1 & 0
\end{pmatrix},\ \
k\mapsto
\begin{pmatrix}
q^{\2} & 0\\ 0 & q^{-\2}
\end{pmatrix}.
$$
Recall next that the $R$-matrix has the form
$$
\RR=q^{2\log_qk\otimes\log_qk}\sum^\infty_{n=0}R_n(q)(ke)^n\otimes
(fk^{-1})^n,
$$
where $\log_q$ is the usual logarithm with base $q$ and $R_0(q)=1$
and $R_1(q)=q-q^{-1}$. Since $(s\,\add_q)(f^n)=0$ for $n\ge2$, we
get
$$
(\iota\otimes s\,\add_q)(\RR)=
\begin{pmatrix}
k & 0\\ q^{-\2}(q-q^{-1})e & k^{-1}
\end{pmatrix},
$$
whence
$$
q^{-\sqrt{2}D_q}=q^{\frac{3}{2}}
\begin{pmatrix}
\partial_q(k^2+q^{-1}(q-q^{-1})^2fe) & q^{-\2}(q-q^{-1})
\partial_q(fk^{-1})\\
 q^{-\2}(q-q^{-1})\partial_q(k^{-1}e) & \partial_q(k^{-2})
\end{pmatrix}.
$$
\end{example}

\bigskip

\section{Spectral triple}

Our next goal is to study commutators of $D_q$ with regular
functions on $G_q$.

\begin{proposition} \label{Commutator1}
For any $a\in\C[G_q]$ we have
$$
[D_q,\pi_{r,q}(a)\otimes1]=-(\pi_{r,q}(a_{(0)})\otimes1)
(\partial_q\varphi^{-1}\otimes s)
(a_{(1)}\varphi^{-1}\otimes\iota\otimes\iota)(UTU^*),
$$
where
$$
U=(\iota\otimes\iota\otimes\add)((\F\otimes1)
(\Dhat\otimes\iota)(\F))
$$
is a unitary operator in $W^*(G)\bar\otimes W^*(G)\otimes\Clg$, the
operator $T\in\U(G\times G)\otimes\Clg$ is defined by
$$
T=(\iota\otimes\iota\otimes\gamma)(t_{13})
+(\iota\otimes\iota\otimes\gamma)(t_{23})
-(\iota\otimes\iota\otimes\add)(\Phi^*)
(\iota\otimes\iota\otimes\gamma)(t_{23})
(\iota\otimes\iota\otimes\add)(\Phi),
$$
and
$$
\Phi=(\iota\otimes\Dhat)(\F^*)
(1\otimes\F^*)(\F\otimes1)(\Dhat\otimes\iota)(\F) \in
W^*(G)\bar\otimes W^*(G)\bar\otimes W^*(G)
$$
is the associator defined by the unitary twist $\F$.
\end{proposition}

Note that in the case $q=1$ we can take $\F=1$, then $U=1$,
$\Phi=1$, $T=(\iota\otimes\iota\otimes\gamma)(t_{13})=-\sum_k
x_k\otimes 1\otimes \gamma(x_k)$, and we recover the familiar
formula
$$
[D,\pi_r(a)\otimes1]=\sum_k
a_{(1)}(x_k)\pi_r(a_{(0)})\otimes s\gamma(x_k)=c(da).
$$

\bp[Proof of Proposition~\ref{Commutator1}] Since
$\pi_{r,q}(a)\otimes1$ commutes with $\sum_k1\otimes
s(\gamma(x_k)\add(x_k))$, it is only the part
$-(\iota\otimes\gamma)(t)$ of $\D$ which contributes to the
commutator. Thus the commutator is the difference of
\begin{equation} \label{ecom2}
(\pi_{r,q}(a)\otimes1)(\partial_q\varphi^{-1}\otimes s)
\Big((\iota\otimes\add)(\F)(\iota\otimes\gamma)(t)
(\iota\otimes\add)(\F^*)\Big)
\end{equation}
and
\begin{equation} \label{ecom1}
(\partial_q\varphi^{-1}\otimes s)\Big((\iota\otimes\add)(\F)
(\iota\otimes\gamma)(t) (\iota\otimes\add)(\F^*)\Big)
(\pi_{r,q}(a)\otimes1).
\end{equation}
Applying (\ref{ecom1}) to a vector $\pi_{r,q}(b)\xi_q\otimes\zeta$
with $b\in\C[G_q]$ and $\zeta\in\Sp$, by definition (\ref{eRegL}) of
$\partial_q$ we get
\begin{equation} \label{ecom3}
(a_{(1)}b_{(1)}\varphi^{-1}\otimes s)\Big((\iota\otimes\add)(\F)
(\iota\otimes\gamma)(t)
(\iota\otimes\add)(\F^*)\Big)(\pi_{r,q}(a_{(0)}b_{(0)})
\xi_q\otimes\zeta)
\end{equation}
For any $c,d\in\C[G_q]$ we have

$\displaystyle
(cd\varphi^{-1}\otimes\iota)\Big((\iota\otimes\add)(\F)
(\iota\otimes\gamma)(t) (\iota\otimes\add)(\F^*)\Big)$
\begin{align*}
&=(c\otimes d\otimes\iota)(\Dhat_q\varphi^{-1}\otimes\iota)
\Big((\iota\otimes\add)(\F)(\iota\otimes\gamma)(t)
(\iota\otimes\add)(\F^*)\Big)\\
&=(c\varphi^{-1}\otimes d\varphi^{-1}\otimes\iota)\Big((\F\otimes 1)
(\Dhat\otimes\iota)
\Big((\iota\otimes\add)(\F)(\iota\otimes\gamma)(t)
(\iota\otimes\add)(\F^*)\Big)(\F^*\otimes1)\Big)\\
&=(c\varphi^{-1}\otimes d\varphi^{-1}\otimes\iota)
(U(\Dhat\otimes\gamma)(t)U^*).
\end{align*}
Therefore (\ref{ecom3}) equals
\begin{multline*}
(a_{(1)}\varphi^{-1}\otimes b_{(1)}\varphi^{-1}\otimes s)
\Big(U(\Dhat\otimes\gamma)(t)U^*\Big)
(\pi_{r,q}(a_{(0)}b_{(0)})\xi_q\otimes\zeta)\\
=(\pi_{r,q}(a_{(0)})\otimes1) (a_{(1)}\varphi^{-1}\otimes
\partial_q\varphi^{-1}\otimes s)
\Big(U(\Dhat\otimes\gamma)(t)U^*\Big)
(\pi_{r,q}(b)\xi_q\otimes\zeta).
\end{multline*}
In other words, (\ref{ecom1}) is equal to
\begin{equation} \label{ecom4}
(\pi_{r,q}(a_{(0)})\otimes1) (\partial_q\varphi^{-1}\otimes s)
(a_{(1)}\varphi^{-1}\otimes\iota\otimes\iota)
(U(\Dhat\otimes\gamma)(t)U^*).
\end{equation}

Consider now the operator (\ref{ecom2}). We can write it as
$$
(\pi_{r,q}(a_{(0)})\otimes1)(\partial_q\varphi^{-1}\otimes s)
(a_{(1)}\varphi^{-1}\otimes\iota\otimes\iota)
\Big((\iota\otimes\iota\otimes\add)(1\otimes\F)
(\iota\otimes\iota\otimes\gamma)(t_{23})
(\iota\otimes\iota\otimes\add)(1\otimes\F^*)\Big).
$$
Since $(\iota\otimes\iota\otimes\gamma)(t_{23})$ commutes with
$(\iota\otimes\iota\otimes\add)(\iota\otimes\Dhat)(\F)$, instead of
conjugating $(\iota\otimes\iota\otimes\gamma)(t_{23})$ by
$(\iota\otimes\iota\otimes\add)(1\otimes\F)$ in the above
expression, we can conjugate it by
$$
(\iota\otimes\iota\otimes\add)((1\otimes\F)(\iota\otimes\Dhat)(\F))
=U(\iota\otimes\iota\otimes\add)(\Phi^*).
$$
Thus (\ref{ecom2}) equals
\begin{equation} \label{ecom5}
(\pi_{r,q}(a_{(0)})\otimes1)(\partial_q\varphi^{-1}\otimes s)
(a_{(1)}\varphi^{-1}\otimes\iota\otimes\iota)
\Big(U(\iota\otimes\iota\otimes\add)(\Phi^*)
(\iota\otimes\iota\otimes\gamma)(t_{23})
(\iota\otimes\iota\otimes\add)(\Phi)U^*\Big).
\end{equation}

To summarize, the commutator $[D_q,\pi_{r,q}(a)\otimes1]$ is equal
to the difference of (\ref{ecom5}) and (\ref{ecom4}). Since
$(\Dhat\otimes\iota)(t)=t_{13}+t_{23}$, this is exactly what the
proposition states. \ep

\begin{corollary} \label{CommAss}
The commutator $[D_q,\pi_{r,q}(a)\otimes1]$ is bounded for all
$a\in\C[G_q]$ if and only if the commutator
$$
[(\pi\otimes\iota\otimes\gamma)(t_{23}),
(\pi\otimes\iota\otimes\add)(\Phi)]
$$
is bounded for any finite dimensional representation $\pi\colon G\to
B(V_\pi)$.
\end{corollary}

One can equivalently formulate the above condition by saying that
the operator
$$
[1\otimes\D,(\pi\otimes\iota\otimes\add)(\Phi)]
$$
affiliated with $B(V_\pi)\otimes W^*(G)\otimes\Clg$ is bounded for
any finite dimensional representation $\pi$ of $G$.

\begin{proof}[Proof of Corollary~\ref{CommAss}]
First observe that $[D_q,\pi_{r,q}(a)\otimes1]$ is bounded for all
$a\in\C[G_q]$ if and only if the operator
$$
(a\varphi^{-1}\otimes\iota\otimes s)(UTU^*)\in \U(G)\otimes B(\Sp)
$$
is bounded for all $a\in\C[G_q]$. Indeed, it is clear that
boundedness of such operators implies boundedness of the
commutators. Conversely, assume that all the commutators are
bounded, and write $1\otimes a$ as a finite sum of elements of the
form $(b\otimes1)\Delta_q(c)$ with $b,c\in\C[G_q]$. Since
$\pi_{r,q}(1)=1$ and
$$
(\pi_{r,q}(bc_{(0)})
\otimes1)(\partial_q\varphi^{-1}\otimes s)
(c_{(1)}\varphi^{-1}\otimes\iota\otimes\iota)(UTU^*)
$$
is bounded by assumption, we conclude that
$$
(a\varphi^{-1}\otimes\partial_q\varphi^{-1}\otimes s)(UTU^*)
$$
is bounded. Then $(a\varphi^{-1}\otimes\iota\otimes s)(UTU^*)$ is
bounded as the representation $\partial_q\varphi^{-1}$ is faithful.

Next notice that when $a$ runs through all elements of $\C[G_q]$,
the functionals $a\varphi^{-1}$ run through the linear span of
matrix coefficients of all finite dimensional representations of
$G$. So to say that $(a\varphi^{-1}\otimes\iota\otimes s)(UTU^*)$ is
bounded for all $a\in\C[G_q]$ is the same as saying that
$$
(\pi\otimes\iota\otimes s)(UTU^*)
$$
is bounded for any finite dimensional unitary representation $\pi$
of $G$. Since $(\pi\otimes\iota\otimes s)(U)$ is unitary, this in
turn is equivalent to boundedness of $(\pi\otimes\iota\otimes
s)(T)$.

Now consider the expression for $(\pi\otimes\iota\otimes s)(T)$. The
first term $(\pi\otimes\iota\otimes s\gamma)(t_{13})$ is clearly
bounded. On the other hand, since $(\pi\otimes\iota\otimes
s\,\add)(\Phi)$ is unitary, the remaining part of
$(\pi\otimes\iota\otimes s)(T)$ can be written as
$$
-(\pi\otimes\iota\otimes s\,\add)(\Phi^*)
[(\pi\otimes\iota\otimes s\gamma)(t_{23}),
(\pi\otimes\iota\otimes s\,\add)(\Phi)].
$$
Therefore the commutators $[D_q,\pi_{r,q}(a)\otimes1]$
are bounded if and only if
$$
(\iota\otimes\iota\otimes s)(
[(\pi\otimes\iota\otimes\gamma)(t_{23}),
(\pi\otimes\iota\otimes \add)(\Phi)])
$$
is bounded for any $\pi$. This is what we need if $\g$ is even
dimensional, as the representation $s$ is then faithful. In the odd
dimensional case there exists another irreducible representation
$\tilde s\colon\Clg\to B(\Sp)$. Then $s\oplus\tilde s$ is faithful.
The representations $s$ and $\tilde s$ are equivalent when
restricted to the even subalgebra of $\Clg$. It follows that $s$ is
isometric on the even subalgebra. But then we conclude that it is
also isometric on the odd part of $\Clg$ by observing that
$\|x\|^2=\|x^*x\|$ and if $x\in\Clg$ is odd then $x^*x$ is even.
Now note that $(\iota\otimes\iota\otimes\gamma)
(t_{23})$ is odd, while $(\pi\otimes\iota\otimes \add)(\Phi)$ is
even, so their commutator is odd. \ep

We now want to get an estimate of the norms of the above commutators
in the case of the Drinfeld associator. But first we establish a
couple of commutation relations.

\begin{lemma} \label{Gaudin1}
We have $[(\iota\otimes\gamma)(t),(\iota\otimes\add)(t)]=0$.
\end{lemma}

\bp By (\ref{eCas1}) we have
$$
2(\iota\otimes\add)(t)=(\iota\otimes\add)\Dhat(C) -C\otimes
1-1\otimes\add(C).
$$
As we know, $(\iota\otimes\gamma)(t)$ commutes with any element of
the form $(\iota\otimes\add)\Dhat(x)$. Thus it commutes with the
first term on the right hand side of the above identity. It also
clearly commutes with the second term. Finally, as we already
remarked prior to Proposition~\ref{Weitz}, the third term is a
scalar. \ep

It is well-known and easy to check that $[t_{ik},t_{ij}+t_{jk}]=0$
for nonequal $i,j,k$. In the spin representation we have a similar
relation.

\begin{lemma} \label{Gaudin2}
We have $[(\iota\otimes\iota\otimes\gamma)(t_{13}),
(\iota\otimes\iota\otimes\add)(t_{12}+t_{23})]=0$.
\end{lemma}

\bp Applying the flip to the first two factors, we can equivalently
check
\begin{equation*}\label{eGaudin1}
[(\iota\otimes\iota\otimes\gamma)(t_{23}),
(\iota\otimes\iota\otimes\add)(t_{12}+t_{13})]=0.
\end{equation*}
This follows immediately from $t_{12}+t_{13}=(\iota\otimes\Dhat)(t)$.
\ep

The relations $[t_{ik},t_{ij}+t_{jk}]=0$ imply consistency of the
Knizhnik-Zamolodchikov equations, or equivalently, mutual
commutativity of the Hamiltonians of the Gaudin model~\cite{G}.
Similarly, using the two previous lemmas we get the following.

\begin{lemma} \label{Gaudin3}
For any $z\in\C$ we have
$$
[(\iota\otimes\iota\otimes\gamma)((1-z)t_{13}+t_{23}),
(\iota\otimes\iota\otimes\add)((z-1)t_{12}+zt_{23})]=0.
$$
\end{lemma}

\bp By Lemma~\ref{Gaudin2} we have
$$
(1-z)[(\iota\otimes\iota\otimes\gamma)(t_{13}),
(\iota\otimes\iota\otimes\add)((z-1)t_{12}+zt_{23})]
=(z-1)[(\iota\otimes\iota\otimes\gamma)(t_{13}),
(\iota\otimes\iota\otimes\add)(t_{12})].
$$
By  Lemma~\ref{Gaudin1} we also have
$$
[(\iota\otimes\iota\otimes\gamma)(t_{23}),
(\iota\otimes\iota\otimes\add)((z-1)t_{12}+zt_{23})]
=(z-1)[(\iota\otimes\iota\otimes\gamma)(t_{23}),
(\iota\otimes\iota\otimes\add)(t_{12})].
$$
So we just have to check that
$$
[(\iota\otimes\iota\otimes\gamma)(t_{13}+t_{23}),
(\iota\otimes\iota\otimes\add)(t_{12})]=0.
$$
This is indeed true as $t_{13}+t_{23}=(\Dhat\otimes\iota)(t)$.
\ep

We can now prove our main technical result.

\begin{proposition} \label{MainEst}
If $\Phi=\Phi_{KZ}$ is the Drinfeld associator, then for any finite
dimensional unitary representation $\pi\colon G\to B(V_\pi)$ we have
$$
\|[(\pi\otimes\iota\otimes\gamma)(t_{23}),
(\pi\otimes\iota\otimes\add)(\Phi)]\| \le
6\|(\pi\otimes\gamma)(t)\|.
$$
\end{proposition}

\bp Fix finite dimensional unitary representations $\pi$ and $\pi'$
of $G$. Put
$$
A=(\pi\otimes\pi'\otimes\add)(t_{12})\ \ \hbox{and}\ \
B=(\pi\otimes\pi'\otimes\add)(t_{23}).
$$
According to (\ref{eKZ}) we have
$$
(\pi\otimes\pi'\otimes\add)(\Phi) =\lim_{a\to 0^+}a^{-\hbar
B}G_a(1-a) a^{\hbar A},
$$
where $G_a$ is such that $G_a(a)=1$ and
$$
G'_a(x)=\hbar\left(\frac{A}{x}+\frac{B}{x-1}\right)G_a(x).
$$
Since all three operators $a^{-\hbar B}$, $G_a(1-a)$ and $a^{\hbar
A}$ are unitary, to prove the proposition is suffices to show that
\begin{equation} \label{ecom10}
[(\pi\otimes\pi'\otimes\gamma)(t_{23}), a^{-\hbar B}]=0,
\end{equation}
\begin{equation}\label{ecom11}
\|[(\pi\otimes\pi'\otimes\gamma)(t_{23}), G_a(x)]\| \le
4\|(\pi\otimes\gamma)(t)\|,
\end{equation}
\begin{equation}\label{ecom12}
\|[(\pi\otimes\pi'\otimes\gamma)(t_{23}), a^{\hbar A}]\| \le
2\|(\pi\otimes\gamma)(t)\|
\end{equation}
for all $a,x\in(0,1)$.

Equality (\ref{ecom10}) follows from Lemma~\ref{Gaudin1}.

To show (\ref{ecom12}) recall that $A$ commutes with
$$
(\pi\otimes\pi'\otimes\gamma)(t_{13}+t_{23})
=(\pi\otimes\pi'\otimes\gamma)(\Dhat\otimes\iota)(t).
$$
Hence the left hand side of (\ref{ecom12}) equals
$\|[(\pi\otimes\pi'\otimes\gamma)(t_{13}), a^{\hbar A}]\|$. Since
$a^{\hbar A}$ is unitary, the latter norm is not larger than
$$
2\|(\pi\otimes\pi'\otimes\gamma)(t_{13})\|=
2\|(\pi\otimes\gamma)(t)\|.
$$

Turning to (\ref{ecom11}), for a fixed $a\in(0,1)$ consider the
commutator
$$
L(x)=[(\pi\otimes\pi'\otimes\gamma)((1-x)t_{13}+t_{23}), G_a(x)].
$$
Then, as $G_a(x)$ is unitary and so
$\|[(\pi\otimes\pi'\otimes\gamma)(t_{13}), G_a(x)]\| \le
2\|(\pi\otimes\gamma)(t)\|$, it is enough to check that
\begin{equation}\label{ecom13}
\|L(x)\|\le 2\|(\pi\otimes\gamma)(t)\|
\end{equation}
for all $x\in(0,1)$. We have for the derivative of $L$ that
$$
L'(x)=\left[(\pi\otimes\pi'\otimes\gamma)((1-x)t_{13}+t_{23}),
\hbar\left(\frac{A}{x}+\frac{B}{x-1}\right)G_a(x) \right]
-[(\pi\otimes\pi'\otimes\gamma)(t_{13}), G_a(x)].
$$
Since $(\pi\otimes\pi'\otimes\gamma)((1-x)t_{13}+t_{23})$ commutes
with $\frac{A}{x}+\frac{B}{x-1}$ by Lemma~\ref{Gaudin3}, we thus see
that $L$ satisfies the differential equation
$$
L'(x)=\hbar\left(\frac{A}{x}+\frac{B}{x-1}\right)L(x)
-[(\pi\otimes\pi'\otimes\gamma)(t_{13}), G_a(x)]
$$
with initial condition $L(a)=0$. Consequently
$$
L(x)=-\int^x_a G_y(x) [(\pi\otimes\pi'\otimes\gamma)(t_{13}),
G_a(y)]dy,
$$
from which we get (\ref{ecom13}) using again unitarity of $G_y(x)$
and $G_a(y)$. \ep

Therefore we get an equivariant spectral triple
$(\C[G_q],L^2(G_q)\otimes\Sp,D_q)$. Since $D_q$ is unitarily
equivalent to $D$, this spectral triple has the same summability
properties as the classical one.

Recall next that if $\g$ is of even dimension $2m$ then the
classical spectral triple is graded by the chirality element
$\chi=i^m\gamma(x_1)\dots\gamma(x_{2m})$. Since $1\otimes\chi$
anticommutes with $\D$ and commutes with elements in the image of
$\iota\otimes\add$, we see that $1\otimes\chi$ anticommutes with
$\D_q$, so our spectral triple for~$G_q$ is even.

To summarize, we have the following result.

\begin{theorem}
If the operator $D_q$ is defined using a unitary Drinfeld twist then
$$
(\C[G_q],L^2(G_q)\otimes\Sp,D_q)
$$
is an equivariant spectral triple of the same parity as the
dimension of $G$.
\end{theorem}

\begin{remark}\mbox{\ }\newline
(i) According to Remark~\ref{Sleb} we can use connections
$\nabla_{f\otimes x}=f\partial(x)\otimes 1+\lambda f\otimes\ad(x)$,
$\lambda\in\R$, to define Dirac operators. By the same procedure as
before we then get Dirac operators on~$G_q$. Since the commutators
of these operators with $\pi_{r,q}(a)\otimes1$, $a\in\C[G_q]$, do
not depend on~$\lambda$, the above theorem remains true for all such
operators.

\noindent (ii) The Dirac operator on $G$ is closely related to Dirac
operators on homogeneous spaces. Consider a homogeneous space $G/K$.
Fix a spin module $\Sp_\kk$ for $\Clk$. Denote by $D_\kk$ the Dirac
operator on the spin bundle over $G/K$ twisted by the bundle induced
by the representation $K\to B(\Sp_\kk)$. On the other hand, we can
can consider the restriction of the Dirac operator on $G$ to the
space of $K$-invariant sections. These two operators can be
expressed in terms of each other. The relation is most transparent
when $G/K$ is even dimensional; then the spaces on which these
operators act can be identified, and the difference of the operators
is bounded.

Therefore using the Dirac operator $D_q$ on $G_q$ we obtain spectral
triples on quantum homogeneous spaces. These are deformations of the
above twisted Dirac operators (whenever they exist). In particular, if $K=T$ is the
maximal torus we obtain a spectral triple on the quantum full flag
manifold $G_q/T$ that is a deformation of the direct sum of
$2^{[\operatorname{rank}\g/2]}$ copies of the spectral triple 
on~$G/T$.
\end{remark}

\end{document}